\documentclass[12pt]{article}
\usepackage{latexsym,amsfonts,amssymb,amsmath}
\usepackage{color}
\usepackage{colordvi,multicol}
\usepackage{graphicx}
\usepackage{float}
\usepackage{subfigure}
\usepackage{epstopdf}
\usepackage{hyperref}
\usepackage{cleveref}
\usepackage{fancyhdr}
\usepackage{lipsum}
\usepackage{tabularx}
\usepackage{caption}
\graphicspath{{figures/}}

\usepackage{setspace}
\newtheorem{thm}{Theorem}[section]

\newtheorem{lem}[thm]{Lemma}

\newtheorem{rem}[thm]{Remark}

\date{}

\begin{document}

\title{\bf A study on the $F$-distribution motivated by Chv\'{a}tal's theorem}

 \author{Qian-Qian Zhou$^{1,}$\footnote{Corresponding author: qianqzhou@yeah.net}\ , Peng Lu$^2$, Ze-Chun Hu$^{2}$\\ \\
 {\small $^1$ School of Science, Nanjing University of Posts and Telecommunications, Nanjing  210023, China}\\
  {\small $^2$ College of Mathematics, Sichuan  University,
 Chengdu 610065, China}}
\maketitle

\makeatletter

\begin{abstract}

Let $X_{d_1, d_2}$ be an $F$-random variable with parameters $d_1$ and $d_2,$ and  expectation $E[X_{d_1, d_2}]$. In this paper,  for any $\kappa>0,$   we investigate  the infimum value of the probability $P(X_{d_1, d_2}\leq \kappa E[X_{d_1, d_2}])$.  Our motivation comes from
Chv\'{a}tal's theorem on the binomial distribution.
\end{abstract}



\section{Introduction and main result}

 Denote by $B(n,p)$ a binomial random variable  characterized by the parameters $n$ and $p.$
 Chv\'{a}tal's theorem states  that for each definite  $n\geq 2$ and for $m=0,\ldots,n$, the probability $P(B(n,m/n)\leq m)$ hits its  minimum value when  $m$ is closest to $\frac{2n}{3}$.
 The development of Chv\'{a}tal's theorem can be traced back to \cite{BPR21},  \cite{Ja21} and \cite{Su21}. With regard to its applications in machine learning, please refer to  \cite{Do18} and  \cite{GM14}.

 Inspired by Chv\'{a}tal's theorem, a series of papers have investigated   the problem of determining the infimum value of
  the probability $P(X\leq E[X])$, where $X$  stands for    a class of  random variables. The distributions of $X$ include the Poisson,  geometric  and  Pascal distributions (\cite{LXH23}), the Gamma distribution (\cite{SHS-a}), the Weibull  and  Pareto distributions (\cite{LHZ23}), the negative binomial distribution (\cite{GTH23}), and some infinitely divisible distributions (\cite{HLZZ24}).

  In  this paper, we continue to study another  basic and important probability distribution: $F$-distribution. Let $d_1\in \mathbb{N},2<d_2\in \mathbb{N}$, and  $X_{d_1, d_2}$ be an $F$-random variable parameterized by $d_1$ and $d_2$. Then by \cite{S73},  we have
\begin{align}\label{2.1}
E[X_{d_1, d_2}]=\frac{d_2}{d_2-2},\   P\left ( X_{d_1, d_2}\leq x \right )= I_{\frac{d_1x}{d_1x+d_2}}(\frac{d_1}{2},\frac{d_2}{2}), \ \ x >0,
\end{align}
where $I$ is  the regularized incomplete beta function defined  in the following way
\begin{align}\label{2.2}
I_x(a, b):=\frac{B(x; a, b)}{B(a, b)}:=\frac{\int_{0}^{x}t^{a-1}(1-t)^{b-1}dt}{B(a, b)},
\end{align}
and $B(a,b)$ is the Beta function, i.e., $B(a,b)=\int_0^1x^{a-1}(1-x)^{b-1}dx$.

The main result  is as follows.
\begin{thm}\label{thm-2.1}
Let $X_{d_1, d_2}$  be an $F$-random variable characterized by  parameters  $d_1\in \mathbb{N}$  and $2<d_2\in \mathbb{N}$. Then
\begin{eqnarray*}
\inf_{d_1, d_2}P\left(X_{d_1, d_2}\le \kappa E[X_{d_1,d_2}]\right)
\left\{
\begin{array}{ll}
=0,& \mbox{if}\quad 0<\kappa<1,\\
=\frac{1}{2},&\mbox{if}\quad  \kappa=1,\\
\geq \frac{1}{2},&\mbox{if}\quad  \kappa>1.
\end{array}
\right.
\end{eqnarray*}
In particular, we have
$$
P(X_{d_1, d_2}\le  E[X_{d_1,d_2}])>\frac{1}{2},
$$
and
$$
\inf_{d_1,d_2}P(X_{d_1, d_2}\le  E[X_{d_1,d_2}])=\frac{1}{2}.
$$
\end{thm}
\begin{rem}
(i) Theorem \ref{thm-2.1} tells us that
$$
P\left(X_{d_1, d_2}\le \kappa E[X_{d_1,d_2}]\right)>\frac{1}{2},\ \forall d_1\in \mathbb{N},\ \forall 2<d_2\in \mathbb{N},\ \forall \kappa\geq 1,
$$
which shows that the $F$-distribution possesses measure concentration in some sense.

(ii) In \cite{SHS-a,SHS-b,SHS-c,ZHS}, the following measure concentration of distributions has been studied:
$$
\inf_{\alpha}P\left(|X_{\alpha}-E[X_{\alpha}]|\leq \sqrt{Var(X_{\alpha})}\right)>0,
$$
where $Var(X_{\alpha})$ is the variance of $X_{\alpha}$.

(iii) In \cite{HST}, Hu et al. investigated the anti-concentration function $(0, \infty)\ni y\to \inf_{\tau}P\{|X_{\tau}-{\rm E}[X_{\tau}]|\geq y \sqrt{{\rm Var}(X_{\tau})}\}$. It was showed that the anti-concentration function is not identically zero for some familiar distributions.
\end{rem}

\section{Proof}\setcounter{equation}{0}

For  $\kappa>0,d_1\in \mathbb{N}$ and $2<d_2\in \mathbb{N}$,   define
\begin{align*}
a:=\frac{d_1}{2},\ b:=\frac{d_2}{2},\ q(a,b,\kappa):=\frac{\kappa a}{\kappa a+b-1}.
\end{align*}
Then
$$
a\geq \frac{1}{2},\ b>1,\  q(a,b,\kappa)=\frac{\kappa d_1}{\kappa d_1+d_2-2}=
\frac{d_1\cdot\frac{\kappa d_2}{d_2-2}}{d_1\cdot\frac{\kappa d_2}{d_2-2}+d_2}.
$$
In virtue of (\ref{2.1}) and (\ref{2.2}), we denote
\begin{gather}\label{2.4}
 I_{q(a, b, \kappa)}(a,b):= I_{\frac{\kappa d_1}{\kappa d_1+d_2-2}} (\frac{d_1}{2},\frac{d_2}{2})=P\left(X_{d_1, d_2}\le \kappa E[X_{d_1,d_2}]\right)= \frac{\int_{0}^{q(a, b, \kappa)}t^{a-1}(1-t)^{b-1}dt}{B(a, b)}.
\end{gather}

\begin{subsection}{Case $\kappa \leq 1$}

Note that, for any function $f(x,y),\  \  x\in \mathcal{X}, y\in \mathcal{Y},$  we have
$$\inf_{(x,y)\in \mathcal{X} \times \mathcal{Y}}f(x,y)=\inf_{x\in \mathcal{X}}\inf_{y\in \mathcal{Y}}f(x,y)=\inf_{y\in \mathcal{Y}}\inf_{x\in \mathcal{X}}f(x,y).$$

For the purpose of   proving     Theorem  \ref{thm-2.1},  the subsequent  lemmas are required.
\begin{lem}\label{lem-2.2-1}
For $x\in (0,1),a>0$ and $b>0$, we have
\begin{align}\label{lem-2.2-a}
I_x(a,b+1)=I_x(a,b)+\frac{x^a(1-x)^b}{bB(a,b)}.
\end{align}
\end{lem}
{\bf Proof.} By the integration by parts, we have
\begin{align}\label{lem-2.2-b}
\int_0^xt^a(1-t)^{b-1}dt&=-\frac{1}{b}\left[x^a(1-x)^b-a\int_0^xt^{a-1}(1-t)^bdt\right].
\end{align}
We also have
\begin{align}\label{lem-2.2-c}
\int_0^xt^a(1-t)^{b-1}dt&=\int_0^xt^{a-1}[1-(1-t)](1-t)^{b-1}dt\nonumber\\
&=\int_0^xt^{a-1}(1-t)^{b-1}dt-\int_0^xt^{a-1}(1-t)^bdt.
\end{align}
By (\ref{lem-2.2-b}) and (\ref{lem-2.2-c}), we get
$$
\int_0^xt^{a-1}(1-t)^{b-1}dt=-\frac{1}{b}x^a(1-x)^b+\frac{a+b}{b}\int_0^xt^{a-1}(1-t)^bdt.
$$
Then together with  $B(a,b+1)=\frac{b}{a+b}B(a,b)$, we derive   that
\begin{align}\label{lem-2.2-d}
\frac{\int_0^xt^{a-1}(1-t)^{b-1}dt}{B(a,b)}&
=-\frac{1}{bB(a,b)}x^a(1-x)^b+\frac{a+b}{bB(a,b)}\int_0^xt^{a-1}(1-t)^bdt\nonumber\\
&=-\frac{1}{bB(a,b)}x^a(1-x)^b+\frac{\int_0^xt^{a-1}(1-t)^bdt}{B(a,b+1)}.
\end{align}
By (\ref{2.2}) and (\ref{lem-2.2-d}), we know that (\ref{lem-2.2-a}) holds. \hfill\fbox

\begin{lem}\label{lem-1.2}
For   $0<\kappa \leq 1$, $a\geq \frac{1}{2}$ and $b>1$,  we have
\begin{gather}\label{a2}
I_{q(a, b+1, \kappa)}(a, b+1) < I_{q(a, b, \kappa)}(a,b).
\end{gather}
\end{lem}
{\bf Proof.} For  any fixed  $a\geq \frac{1}{2}$,
by  $I_x(a,b+1)=I_x(a,b)+\frac {x^a(1-x)^b}{bB(a,b)}$ and (\ref{2.4}), we have
\begin{align*}
&I_{q(a,b+1, \kappa)}(a,b+1)-I_{q(a,b, \kappa)}(a,b)\\
&= \left(\frac{\int_0^{q(a, b+1, \kappa)}t^{a-1}(1-t)^{b-1}dt}{B(a,b)}+\frac{q^a(a, b+1, \kappa)(1-q(a, b+1, \kappa))^b}{bB(a,b)}\right)\\
 &\quad-\frac{\int_0^{q(a, b, \kappa)}t^{a-1}(1-t)^{b-1}dt}{B(a,b)}\\
&=\frac {1}{bB(a,b)}\left[-b\int_{q(a,b+1, \kappa)}^{q(a,b, \kappa)}t^{a-1}(1-t)^{b-1}dt+{q^a(a,b+1, \kappa)}(1-{q(a,b+1, \kappa)})^b\right].
\end{align*}

Define a function
$$
g(t):=t^a (1-t)^{b-1},\  \  t\in (0,1).
$$
If $g'(t)=0$,  then $t=\frac{a}{a+b-1}\in (0,1)$.  It is simple to acquire   that $g'(t)>0$  on $(0, \frac{a}{a+b-1})$ and $g'(t)<0$  on $(\frac{a}{a+b-1}, 1)$.
  Since  $q(a,b,\kappa)\le q(a,b,1)=\frac{a}{a+b-1}$ for  $0<\kappa \le 1,$     then the function $g(t)$ is strictly  increasing on  $(q(a, b+1, \kappa), q(a,b,\kappa)).$ Therefore, we have
\begin{align}\label{lem2.3-b}
\hspace{-1cm}&\eqref{a2}  \text{  holds }\nonumber\\
 &\Leftrightarrow b\int_{q(a,b+1,\kappa)}^{q(a,b,\kappa)}t^{a-1}(1-t)^{b-1}dt>{q^a(a,b+1,\kappa)}
 \big(1-{q(a,b+1,\kappa)}\big)^b\nonumber\\
&\Leftrightarrow b\int_{q(a,b+1,\kappa)}^{q(a,b,\kappa)}\frac{g(t)}{t}dt>{q^a(a,b+1,\kappa)}
\big(1-{q(a,b+1,\kappa)}\big)^b\nonumber\\
&\Leftarrow  b\int_{q(a,b+1,\kappa)}^{q(a,b,\kappa)}g(t)dt\geq q(a,b,\kappa){q^a(a,b+1,\kappa)}\big(1-{q(a,b+1,\kappa)}\big)^b\nonumber\\
&\Leftarrow  b\big(q(a,b,\kappa)-q(a,b+1,\kappa)\big)g(q(a,b+1,\kappa))\geq q(a,b,\kappa){q^a(a,b+1,\kappa)}\big(1-{q(a,b+1,\kappa)}\big)^b.
\end{align}
According to the definitions  of $q(a, b, \kappa)$  and function $g(\cdot)$,  it follows  that
\begin{align*}
&b\big(q(a,b,\kappa)-q(a,b+1,\kappa)\big)g(q(a,b+1,\kappa))\\
&= \frac{\kappa ab}{(\kappa a+b-1)(\kappa a+b)}\big(\frac{\kappa a}{\kappa a+b}\big)^a\big(\frac{b}{\kappa a +b}\big)^{b-1}\\
&= \frac{\kappa a }{\kappa a+b-1}\big(\frac{\kappa a}{\kappa a+b}\big)^a\big(\frac{b}{\kappa a +b}\big)^{b}\\
&=q(a,b,\kappa){q^a(a,b+1,\kappa)}\big(1-{q(a,b+1,\kappa)}\big)^b.
\end{align*}
Then it  implies that the last inequality in (\ref{lem2.3-b}) holds and thus \eqref{a2}  is true. \hfill\fbox

\begin{lem}\label{lem-1.3}
For  $0<\kappa\le1$  and   $a=\frac{d_1}{2}, d_1\in \mathbb{N}$,  we have
\begin{gather}\label{a4}
g_\kappa (a):=\inf_{b=d_2/2,2<d_2\in \mathbb{N}}  I_{q(a,b, \kappa)}(a,b)=  \lim_{b\to \infty}I_{q(a,b, \kappa)}(a,b)=\frac {\int_{0}^{\kappa a}t^{a-1}e^{-t}dt}{\Gamma (a)}.
\end{gather}
\end{lem}
{\bf Proof.} To prove (\ref{a4}), by Lemma \ref{lem-1.2}, it suffices to demonstrate that
\begin{gather}\label{a4-1}
\lim_{b\to \infty}I_{q(a,b, \kappa)}(a,b)=\frac {\int_{0}^{\kappa a}t^{a-1}e^{-t}dt}{\Gamma (a)}.
\end{gather}

By  $I_{q(a,b,\kappa)}(a,b)$ given  in \eqref{2.4} and the change of variable $t
=\frac {\kappa a}{\kappa a+b-1}y $,  we have
 \begin{eqnarray*}
\lim_{b\to \infty}I_{q(a,b, \kappa)}(a,b)
&=&\lim_{b\to \infty}\frac{\int_0^{\frac {\kappa a}{\kappa a+b-1}}t^{a-1}(1-t)^{b-1}dt}{B(a,b)}\\
&=&\lim_{b\to \infty}\frac{({\frac {\kappa a }{\kappa a+b-1})}^{a}\int_{0}^{1}y^{a-1}\left(1-\frac {\kappa a}{\kappa a+b-1}y\right)^{b-1}dy}{B(a,b)}\\
&=&\lim_{b\to \infty}\frac{\Gamma(a+b)}{\Gamma(a)\Gamma(b)}\left(\frac {\kappa a }{\kappa a+b-1}\right)^{a}\int_{0}^{1}y^{a-1}\left(1-\frac {\kappa a}{\kappa a+b-1}y\right)^{b-1}dy\\
&=&\lim_{b\to \infty}\frac{\Gamma(a+b)}{\Gamma(b)b^a}\left(\frac{\kappa  ab}{\kappa a+b-1}\right)^a\frac{1}{\Gamma(a)}\int_{0}^{1}y^{a-1}\left(1-\frac {\kappa a}{\kappa a+b-1}y\right)^{b-1}dy,
\end{eqnarray*}
where  we used the fact that $B(a,b)=\frac{\Gamma(a)\Gamma(b)}{\Gamma(a+b)}.$

By the special limit $\lim_{x\rightarrow 0}(1+x)^{\frac{1}{x}}=e,$  we have $\lim_{b\rightarrow \infty}(1-\frac {\kappa a}{\kappa a+b-1}y)^{b-1}=e^{-\kappa a y}.$ By $ \Gamma(x+1) \approx \sqrt{2\pi}e^{-x}x^{x+\frac{1}{2}}$ as $x \to \infty,$  we get that $\lim_{b\rightarrow \infty}\frac{\Gamma(a+b)}{\Gamma(b)b^a}=1.$  Then, by the bounded convergence theorem, we derive  that
\begin{eqnarray}\label{d1}
\lim_{b\to \infty}I_{q(a,b, \kappa)}(a,b)&=&\frac{(\kappa a)^a}{\Gamma(a)}\int^1_{0} y^{a-1}e^{-\kappa a y}dy.
\end{eqnarray}
By letting $t=\kappa a y$ in \eqref{d1}, we get
$$\lim_{b\to \infty}I_{q(a,b, \kappa)}(a,b)=\frac {\int_{0}^{\kappa a}t^{a-1}e^{-t}dt}{\Gamma (a)}.$$
Hence (\ref{a4-1}) holds. \hfill\fbox

\bigskip

By Sun et al. \cite[Section 2.1]{SHS-a}, we know that
\begin{eqnarray}\label{2.13}
\inf_{a>0} \frac {\int_{0}^{\kappa a}t^{a-1}e^{-t}dt}{\Gamma (a)}=\left\{
\begin{array}{lc}
0,& \kappa\in (0,1),\\
\frac{1}{2},& \kappa=1.
\end{array}
\right.
\end{eqnarray}
By Lemmas \ref{lem-1.2},  \ref{lem-1.3}, and (\ref{2.13}), we obtain that Theorem \ref{thm-2.1} holds when  $0<\kappa\le 1$.

\end{subsection}

%\newpage

\begin{subsection}{Case $ \kappa > 1 $}
In this subsection, we assume that $\kappa>1.$
%This case is complicated  and interesting.
For any given numbers $a\in \mathcal{A}=\big\{  \frac{1}{2}d_1| d_1 \in \mathbb{N}  \big\}$  and $b\in \mathcal{B}=\big\{  \frac{1}{2}d_2| d_2 \in \mathbb{N}, d_2>2  \big\}$,
it is evident  that the function $q(a, b, \kappa)$  is  strictly increasing with respect to $\kappa.$  Then for any $\kappa>1$, we obtain that
\begin{gather*}
I_{q(a,b, \kappa)}(a,b):= P(X_{d_1, d_2} \leq \kappa E[X_{d_1, d_2}]) >I_{q(a,b, 1)}(a,b):=P(X_{d_1, d_2} \leq  E[X_{d_1, d_2}]),
\end{gather*}
which,  along  with  Section 2.1,  implies  that in this case,
$$
\inf_{d_1, d_2}P\left(X_{d_1, d_2}\le \kappa E[X_{d_1,d_2}]\right)\geq \frac{1}{2}.
$$

By virtue of  Python, we obtain the following numerical results for different values of $\kappa$.

\begin{table}[htbp]
\centering
\captionof{table}{The $\inf P(X_{d_1. d_2} \leq \kappa E[X_{d_1, d_2}]) $ for different values of $\kappa$}
\begin{tabularx}{\textwidth}{X<{\centering}|X<{\centering}|X<{\centering}|X< {\centering}}
  \hline
  $\kappa$ & $\inf P(X_{d_1, d_2} \leq \kappa E[X_{d_1, d_2}])$ & $d_1$ & $d_2$ \\
  \hline
  %0.8 & 3.51912e-07 & 1999 & 1999 \\
  %0.95 & 0.130488 & 1999 & 1999 \\
  %1.0 & 0.508925  & 1999 & 1999 \\
  1.00005 & 0.509371 & 1999 & 1999 \\
  1.001 & 0.516817 & 667 & 1999 \\
  1.005 & 0.533577 & 134 & 1999 \\
  1.05 & 0.601371 & 14 & 1999 \\
  1.5 & 0.776954 & 2 & 1999 \\
  3.0 & 0.936000 & 1 & 1999 \\
  3.005 & 0.916991 & 1 & 803 \\
  3.05 & 0.919240 & 1 & 83 \\
  $\pi$ & 0.923510 & 1 & 31 \\
  4.0 & 0.950133 & 1 & 7 \\
  6.0 & 0.974279 & 1 & 4 \\
  8.0 & 0.983723 & 1 & 3 \\
  16.0 & 0.993835 & 1  & 3 \\
  \hline
\end{tabularx}
\end{table}

%\newpage

%Below are graphs of the function   $P(X_{d_1, d_2}\leq \kappa E[X_{d_1, d_2}])$  for different values of  $\kappa.$
%\vspace { 1cm }

\begin{rem}
Based on the above numerical results, we conjecture that for $\kappa>1$,
$$
\inf_{d_1, d_2}P\left(X_{d_1, d_2}\le \kappa E[X_{d_1,d_2}]\right)> \frac{1}{2}.
$$
Up to now, we can not give the proof.
\end{rem}

\end{subsection}

\noindent {\bf\large Acknowledgments}\quad  This work was supported by NNSFC (12301603, 12171335).


\begin{thebibliography}{1234}

\bibitem{BPR21} L. Bababesi, L. Pratelli, P. Rigo, On the Chv\'{a}tal-Janson conjecture. Statis. Probab. Lett., 2023,  194:  109744.

\bibitem{Do18} B. Doerr,  An elementary analysis of the probability that a binomial random variable exceeds its expectation. Statis. Probab. Lett., 2018, 139:  67-74.

\bibitem{GM14} S. Greenberg, M. Mohri,  Tight lower bound on the probability of a binomial exceeding its expectation.   Statis. Probab. Lett., 2014, 86: 91-98.

\bibitem{GTH23}  Z.-Y. Guo, Z.-Y.  Tao, Z.-C. Hu,  A study on the negative binomial distribution motivated by Chv\'{a}tal's theorem. Statis. Probab. Lett.,  2024, 207: 110037.

\bibitem{HLZZ24} Z.-C. Hu, P. Lu, Q.-Q. Zhou, X.-W. Zhou, The infimum values of the probability functions for some infinitely divisible distributions  motivated by Chv\'{a}tal's theorem. J. of Math. (PRC), 2024, 44(4): 309-316.

\bibitem{HST} Z.-C. Hu, R. Song, Y. Tan, On the anti-concentration functions of some familiar families of distributions.   Math. Theory Appl., 2024, 44: 1-15.

\bibitem{Ja21} S. Janson, On the probability that a binomial variable is at most its expectation.  Statis. Probab. Lett., 2021, 171: 109020.

\bibitem{LHZ23} C. Li, Z.-C. Hu,  Q.-Q. Zhou,  A study on the Weibull and Pareto distributions motivated by Chv\'{a}tal's  theorem.  J. of Math. (PRC), 2024, 44(3): 195-202.

\bibitem{LXH23} F.-B. Li, K. Xu, Z.-C. Hu, A study on the Poisson, geometric and Pascal distributions motivated by Chv\'{a}tal's conjecture. Statis. Probab. Lett., 2023, 200: 109871.

%\bibitem{PR16} C. Pelekis, J. Ramon,  A lower bound on the probability that a binomial random variable is exceeding its mean.  Statis. Probab. Lett., 2016,  119: 305-309.

\bibitem{S73}   F.-W. Steutel,  Some recent results in infinite divisibility. Stoch. Proc. Appl., 1973, 1: 125-143.

\bibitem{Su21} P. Sun,  Strictly unimodality of the probability that the binomial distribution is more than its expectation. Discrete Appl. Math., 2021, 301:  1-5.

\bibitem{SHS-a} P. Sun, Z.-C. Hu, W. Sun, The infimum  values of two probability functions for the Gamma distribution. J. Inequal. Appl.,  2024, 2024(1): 5.

\bibitem{SHS-b} P. Sun, Z.-C. Hu, W. Sun, Variation comparison between infinitely divisible distributions and the normal distribution, Accepted by Statistical Papers, 2024.

\bibitem{SHS-c} P. Sun, Z.-C. Hu, W. Sun, Variation comparison between the $F$-distribution and the normal distribution, arXiv: 2305.13615, 2023.

\bibitem{ZHS} J. Zhang, Z.-C. Hu, W. Sun, On the measure concentration of infinitely divisible distributions, Accepted by Acta Mathematica Scientia, 2024.


\end{thebibliography}
\end{document}